# Correcting Newton–Côtes integrals by Lévy areas

IVAN NOURDIN[1] and THOMAS SIMON[2]

[1] *Laboratoire de Probabilités et Modèles Aléatoires, Université de Paris, 4 Place Jussieu, Boîte courrier 188, 75252 Paris Cedex 5, France. E-mail: nourdin@ccr.jussieu.fr*
[2] *Equipe d'Analyse et Probabilités, Université d'Evry-Val d'Essonne, Boulevard François Mitterand, 91025 Evry Cedex, France. E-mail: tsimon@univ-evry.fr*

In this note we introduce the notion of Newton–Côtes functionals corrected by Lévy areas, which enables us to consider integrals of the type $\int f(y)\,\mathrm{d}x$, where $f$ is a $\mathscr{C}^{2m}$ function and $x,y$ are real Hölderian functions with index $\alpha > 1/(2m+1)$ for all $m \in \mathbb{N}^*$. We show that this concept extends the Newton–Côtes functional introduced in Gradinaru *et al.*, to a larger class of integrands. Then we give a theorem of existence and uniqueness for differential equations driven by $x$, interpreted using the symmetric Russo–Vallois integral.

## 1. Introduction

In stochastic modelling, differential equations driven by a fractional Brownian motion $B^H$,

$$X_t = x_0 + \int_0^t b(X_s)\,\mathrm{d}s + \int_0^t \sigma(X_s)\,\mathrm{d}B_s^H, \qquad t \in [0,1], \tag{1.1}$$

are popular generalizations of classical stochastic differential equations (SDEs) driven by Brownian motion, relevant especially in finance; see, for example, Cheridito [2], Cutland, Kopp and Willinger [7], Comte and Renault [4] and the references therein. Since fractional Brownian motion (fBm for short) $B^H$ is a semimartingale if and only if its Hurst index $H$ equals $1/2$ (i.e., when $B^H$ is the standard Brownian motion), for $H \neq 1/2$ the meaning of $\int_0^t \sigma(X_s)\,\mathrm{d}B_s^H$ in (1.1) is not the usual one and has to be made precise. Let us give a brief sketch of the three theories of integration with respect to fBm that are frequently used nowadays, referring to the survey by Coutin [5] for a further account:

(a) Russo and Vallois [19] introduced a symmetric integral defined by

$$\int_0^t Z_s \,\mathrm{d}^\circ B_s^H = \lim_{\varepsilon \to 0} \varepsilon^{-1} \int_0^t \frac{Z_{s+\varepsilon} + Z_s}{2} (B_{s+\varepsilon}^H - B_s^H) \,\mathrm{d}s, \tag{1.2}$$

provided the (uniform in probability) limit exists. When the integrand $Z$ is of the type $Z_s = f(B_s^H)$, Cheridito and Nualart [3] and Gradinaru *et al.* [10] showed that $\int_0^t f(B_s^H) \,\mathrm{d}^\circ B_s^H$ exists for all $f : \mathbb{R} \to \mathbb{R}$ regular enough if and only if $H > 1/6$. When $Z_s = h(B_s^H, V_s)$ with $V$ a process of bounded variation and $h : \mathbb{R}^2 \to \mathbb{R}$ a regular function, it is easily shown in using a Taylor expansion that $\int_0^t h(B_s^H, V_s) \,\mathrm{d}^\circ B_s^H$ exists if $H > 1/3$. When $H \leq 1/3$, one can extend the definition (1.2) and give a meaning to

$$\int_0^t h(B_s^H, V_s) \,\mathrm{d}B_s^H \tag{1.3}$$

with the help of the $m$-order Newton–Côtes functional introduced in Gradinaru *et al.* [10]; see Definition 2.1 thereafter. Choosing $m$ sufficiently large exhibits a stochastic functional that makes sense of (1.1) for *all* $H \in (0,1)$; see Nourdin [14]. However, one needs to suppose somewhat arbitrarily that the solution to (1.1) is a priori of the type $h(B_s^H, V_s)$.

(b) Another formalism relies on the Malliavin calculus for fBm in the sense of Nualart and Zakai (see Nualart [16]), and more specifically on Skorohod's integration operator $\delta^H$. Combining this formalism with fractional calculus techniques and Young integrals, one can then study (1.1) for $H > 1/2$ in any dimension. We refer to Nualart and Rasçanu [18] and the survey article by Nualart [17] for further topics on this theory.

(c) Finally, one can define (1.1) with the help of the rough paths theory pioneered by Lyons [11]. Roughly speaking, the goal of this theory is to give a meaning to quantities such as $\int_\gamma \omega$, where $\omega$ is a differential 1-form and $\gamma$ is a curve that has only Hölder continuous regularity. To do so, it is then necessary to reinterpret (1.1) using a differential 1-form, through the formulation

$$X_t = x_0 + \int_{\gamma([0,t])} \omega \tag{1.4}$$

with $\gamma_t = (B_t^H, t, X_t) \in \mathbb{R}^3$ and $\omega = \sigma(x_3)\,\mathrm{d}x_1 + b(x_3)\,\mathrm{d}x_2$. Recent results by Coutin and Qian [6] and Feyel and De la Pradelle [9] showed that one can solve (1.4) when $H > 1/4$, in any dimension. Rough paths theory has rich ramifications—see the monograph by Lyons and Qian [12]—but requires a formalism that is sometimes heavy.

It is quite natural to ask whether and how these different theories may intertwine. In Alòs and Nualart [1], the following link is established between (a) and (b): fixing a time horizon $T$ and $H \geq 1/2$, if $u$ is a stochastic process that is regular enough (in the sense of Malliavin calculus), then its symmetric integral along $B^H$ exists and is given by

$$\int_0^T u_t \,\mathrm{d}^\circ B_t^H = \delta^H(u) + H(2H-1) \int_0^T \int_0^T D_s^H u_t |t-s|^{2H-2} \,\mathrm{d}s\,\mathrm{d}t,$$

where $D_s^H$ stands stands for the Malliavin derivative and $\delta^H$ stands for the Skorohod integral. In the case $H = 1/2$, one retrieves the classical formula connecting Itô and Stratonovich integrals.

The present note wishes to link (a) and (c). We propose a correction of the Newton–Côtes integrator $\mathrm{d}^{\mathrm{NC},m}$ by some *Lévy areas*, which are the central objects in rough paths theory. Our new integrator $\mathrm{d}^{A,m}$ gives a meaning to

$$\int_0^t f(y_s)\,\mathrm{d}^{A,m}x_s$$

for all $m \in \mathbb{N}^*$ when $f : \mathbb{R} \to \mathbb{R}$ is $\mathscr{C}^{2m}$ and $x, y$ are any fractal functions of index $\alpha > 1/(2m+1)$ (Theorem 2.5). Note that, compared to the works mentioned in (a), our class of integrands is more satisfactory because $y$ does not depend on $x$. Compared to (b) and (c), we also reach a lower level for $H$. However, our approach has the main drawback of being genuinely one dimensional.

In the second part of the paper we focus on the specific case $\alpha > 1/3$. Considering $\mathrm{d}^{A,1}$, which is our generalization of the symmetric Russo–Vallois integrator $\mathrm{d}^\circ$, we prove existence and uniqueness for (1.1) under some standard conditions on the coefficients (Theorem 3.2). The proof relies on Banach's fixed point theorem. Finally, we notice that, for $y_t = g(x_t, \ell_t)$ with $\ell$ of bounded variation, one can choose a Lévy area $A$ of order 0 such that the operators $\mathrm{d}^{A,1}$ and $\mathrm{d}^\circ$ actually coincide (Proposition 3.5). We are not sure whether an identification with Newton–Côtes functionals can be pursued for $m \geq 2$, because of the crucial Chasles relationship in the definition of Lévy areas.

This paper was mainly inspired by Feyel and De la Pradelle [9], more precisely by their first draft. For example, our Lemma 2.7, which is key in establishing Theorem 2.5, can be viewed as a continuous analogue to the sewing lemma, Lemma 2.1, therein. The possibility of reaching any value of $H$ after considering families of Lévy areas was also strongly suggested by Feyel and De la Pradelle [9]. However, our framework is continuous and, in particular, our integrals are true integrals for $H > 1/3$, which may look more natural. Above all, we feel that this formalism is one of the simplest possible, and provides a handy framework for a more advanced analysis of (1.1), examples of which can be found in Neuenkirch and Nourdin [13] and Nourdin and Simon [15].

## 2. Newton–Côtes integrals corrected by Lévy areas

Without loss of generality we will consider functions defined on the interval $[0,1]$. We fix once and for all $m \in \mathbb{N}^*$ and $\alpha \in (1/(2m+1), 1)$. Denote by $\mathrm{C}^\alpha$ the set of fractal functions $z : [0,1] \to \mathbb{R}$ of index $\alpha$, that is, for which there exists $L > 0$ such that for all $s, t \in [0,1]$,

$$|z_t - z_s| \leq L|t-s|^\alpha.$$

Define a family of interpolation measures $\{\nu_m, m \geq 1\}$ by

$$\nu_1 = \tfrac{1}{2}(\delta_0 + \delta_1), \qquad \text{if } m = 1,$$

$$\nu_m = \sum_{j=0}^{2m-2} \left( \int_0^1 \left( \prod_{k\neq j} \frac{2(m-1)u-k}{j-k} \right) \mathrm{d}u \right) \delta_{j/(2m-2)}, \qquad \text{if } m \geq 2,$$

where $\delta$ stands for the Dirac mass. The signed measure $\nu_m$ is the unique discrete measure carried by the numbers $j/(2m-2)$ that coincides with Lebesgue measure when integrated on polynomials of degree smaller than $2m-1$. In Gradinaru *et al.* [10], the Newton–Côtes functional was defined as follows:

***Definition 2.1.*** *Let $x:[0,1]\to\mathbb{R}$, $z:[0,1]\to\mathbb{R}^2$ and $h:\mathbb{R}^2\to\mathbb{R}$ be continuous functions. The quantity defined by*

$$\int_0^t h(z_s)\,\mathrm{d}^{\mathrm{NC},m}x_s \stackrel{\mathrm{def}}{=} \lim_{\varepsilon\to 0} \varepsilon^{-1} \int_0^t \mathrm{d}s\,(x_{s+\varepsilon}-x_s) \int_0^1 h((1-\alpha)z_s+\alpha z_{s+\varepsilon})\nu_m(\mathrm{d}\alpha), \quad (2.5)$$

*provided the limit exists, is called the $m$-order Newton–Côtes functional $I_m(h,z,x)$ of $h(z)$ with respect to $x$.*

***Remarks 2.2.*** (a) When $m=1$, the Newton–Côtes functional is an integral that coincides with the symmetric integral $\int_0^t h(z_s)\,\mathrm{d}^\circ x_s$ given in definition (1.2).

(b) When $m\geq 2$, the Newton–Côtes functional is not an integral, because if $h(z)=\tilde{h}(\tilde{z})$, the identification

$$\int_0^T h(z_s)\,\mathrm{d}^{\mathrm{NC},m}x_s = \int_0^T \tilde{h}(\tilde{z}_s)\,\mathrm{d}^{\mathrm{NC},m}x_s$$

does not hold in general. This explains why we chose the terminology *functional* instead of *integral*.

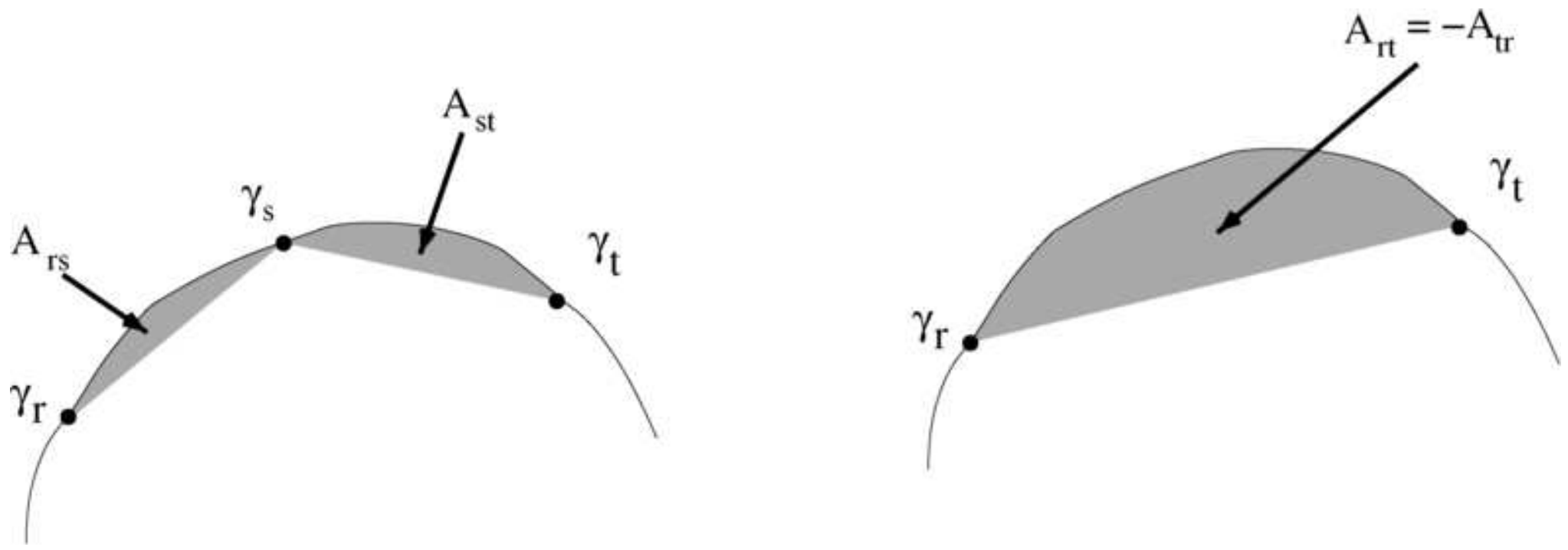

**Figure 1.** Algebraic relation satisfied by a Lévy area of order 0 associated with $\gamma$.

Notice that there is no reason a priori that the functional $I_m(h,z,x)$ exists. Nourdin [14] established existence when $z$ is of the form $u\mapsto f(x_u,\ell_u)$, where $\ell:[0,1]\to\mathbb{R}$ has bounded variations and $f:\mathbb{R}^2\to\mathbb{R}$ is regular enough.

To extend the class of integrands, we wish to define a new concept of functional. To do so, we first introduce the notion of Lévy area. If $x:[0,1]\to\mathbb{R}$ and $y:[0,1]\to\mathbb{R}$ are smooth functions, and if $\gamma:[0,1]\to\mathbb{R}^2$ denotes the curve $\gamma_t=(x_t,y_t)$, the domain $D_{st}$ between the path $\gamma([s,t])$ and the affine chord from $\gamma_s$ to $\gamma_t$ is well defined and we easily verify (see Figure 1) that

$$A_{rs}+A_{st}+A_{tr}=-\mathrm{area}(T_{rst})$$

for all $r,s,t\in[0,1]$, where $A_{st}$ is the (algebraic) area of $D_{st}$ and $T_{rst}$ is the oriented triangle with vertices $\gamma_r$, $\gamma_s$ and $\gamma_t$. Therefore, the following definition is quite natural:

***Definition 2.3.*** *Let $m\in\mathbb{N}^*$, $\alpha\in(0,1)$ and $x,y:[0,1]\to\mathbb{R}$ be two functions in $\mathrm{C}^\alpha$. We say that $A$ is* an $\alpha$-Lévy area of order $2m-2$ associated with $\gamma$ *if for all $s,t\in[0,1]$ the map $P\to A_{st}(P)$ is a linear map from $\mathscr{P}_{2m-2}$ (the space of polynomials in $y$ with degree $\le 2m-2$) into $\mathbb{R}$, if for all $r,s,t\in[0,1],k\in\{0,\dots,2m-2\}$,*

$$A_{rs}(Y^k)+A_{st}(Y^k)+A_{tr}(Y^k)=-\iint_{T_{rst}}\eta^k\,\mathrm{d}\xi\,\mathrm{d}\eta, \tag{2.6}$$

*and if there exists a constant $C>0$ such that for all $s,t\in[0,1],k\in\{0,\dots,2m-2\},\zeta\in[y_s,y_t]$,*

$$|A_{st}[(Y-\zeta)^k]|\le C|t-s|^{2m\alpha}. \tag{2.7}$$

***Remarks 2.4.*** (a) From (2.7), we see that $A_{ss}(P)=0$ for all $s\in[0,1]$ and $P\in\mathscr{P}_{2m-2}$. From (2.6) and because $\iint_{T_{sst}}\eta^k\,\mathrm{d}\xi\,\mathrm{d}\eta=0$, we see that $A_{st}(P)=-A_{ts}(P)$ for all $s,t\in[0,1]$ and $P\in\mathscr{P}_{2m-2}$.

(b) When $k=0$, we write $A$ instead of $A(1)$, and conditions (2.6)–(2.7) become

$$A_{rs}+A_{st}+A_{tr}=-\mathrm{area}(T_{rst}), \tag{2.8}$$

$$|A_{st}|\le c|t-s|^{2\alpha}. \tag{2.9}$$

We can now give the main result and the central definition of this paper:

**Theorem 2.5.** *Let $m\in\mathbb{N}^*$, $x,y\in\mathrm{C}^\alpha$ with $\alpha>1/(2m+1)$ and $A$ be an $\alpha$-Lévy area of order $2m-2$ associated with $\gamma=(x,y)$. For $f:\mathbb{R}\to\mathbb{R}$ a $\mathscr{C}^{2m}$ function, define*

$$I_\varepsilon^\gamma(f)=\varepsilon^{-1}\int_0^1\mathrm{d}u\,(x_{u+\varepsilon}-x_u)\int_0^1 f((1-\alpha)y_u+\alpha y_{u+\varepsilon})\nu_m(\mathrm{d}\alpha)$$
$$+\varepsilon^{-1}\sum_{k=0}^{2m-2}\frac{1}{(k+1)!}\int_0^1 f^{(k+1)}(y_u)A_{u,u+\varepsilon}[(y-y_u)^k]\,\mathrm{d}u$$

*for every* $\varepsilon > 0$. *Then the family* $\{I_\varepsilon^\gamma(f), \varepsilon > 0\}$ *converges when* $\varepsilon \downarrow 0$. *Its limit is denoted by*

$$I^\gamma(f) = \int_0^1 f(y_u)\,\mathrm{d}^{A,m} x_u$$

*and is called the* $m$*-order Newton–Côtes functional corrected by* $A$ *of* $f(y)$ *with respect to* $x$. *In addition, these functionals are compatible in the sense that if* $\alpha > 1/(2m+1)$ *for all* $n > m$ *and any* $\alpha$*-Lévy area* $A$ *of order* $2n-2$ *associated with* $\gamma$, *then*

$$\int_0^1 f(y_u)\,\mathrm{d}^{A,n} x_u = \int_0^1 f(y_u)\,\mathrm{d}^{A,n-1} x_u = \cdots = \int_0^1 f(y_u)\,\mathrm{d}^{A,m} x_u.$$

The proof of Theorem 2.5 relies on the two following simple lemmas.

**Lemma 2.6.** *Let* $m \in \mathbb{N}^*$, $x, y \in \mathrm{C}^\alpha$ *with* $\alpha > 1/(2m+1)$, $A$ *be an* $\alpha$*-Lévy area of order* $2m-2$ *associated with* $\gamma = (x, y)$ *and* $f : \mathbb{R} \to \mathbb{R}$ *be a* $\mathscr{C}^{2m}$ *function. Set*

$$\begin{aligned} I_n(\varepsilon) &= 2^n \varepsilon^{-1} \int_0^{\varepsilon[1/\varepsilon]} \mathrm{d}u\,(x_{u+\varepsilon 2^{-n}} - x_u) \int_0^1 f((1-\alpha) y_u + \alpha y_{u+\varepsilon 2^{-n}}) \nu_m(\mathrm{d}\alpha) \\ &\quad + 2^n \varepsilon^{-1} \sum_{k=0}^{2m-2} \frac{1}{(k+1)!} \int_0^{\varepsilon[1/\varepsilon]} f^{(k+1)}(y_u) A_{u,u+\varepsilon 2^{-n}}[(y - y_u)^k]\,\mathrm{d}u \end{aligned}$$

*for every* $\varepsilon > 0$ *and* $n \in \mathbb{N}$. *The sequence of functions* $\{I_n, n \in \mathbb{N}\}$ *converges uniformly on each compact of* $]0, 1]$, *and the limit* $I_\infty$ *satisfies*

$$I_\infty(\varepsilon) = I_\varepsilon^\gamma(f) + O(\varepsilon^{[(2m+1)\alpha - 1] \wedge \alpha}). \tag{2.10}$$

**Proof.** First, assume that $m = 1$. In this case, we have

$$I_n(\varepsilon) = 2^n \varepsilon^{-1} \left( \int_0^{\varepsilon[1/\varepsilon]} \frac{f(y_u) + f(y_{u+\varepsilon 2^{-n}})}{2} (x_{u+\varepsilon 2^{-n}} - x_u)\,\mathrm{d}u + \int_0^{\varepsilon[1/\varepsilon]} f'(y_u) A_{u,u+\varepsilon 2^{-n}}\,\mathrm{d}u \right),$$

where for simplicity we have written $A_{st}$ instead of $A_{st}(1)$. Decomposing the integral into dyadic intervals and making a change of variable, we first get

$$\begin{aligned} I_n(\varepsilon) &= 2^n \varepsilon^{-1} \sum_{k=0}^{[1/\varepsilon]2^n - 1} \int_{\varepsilon k 2^{-n}}^{\varepsilon(k+1)2^{-n}} \left[ \frac{f(y_u) + f(y_{u+\varepsilon 2^{-n}})}{2} (x_{u+\varepsilon 2^{-n}} - x_u) \right. \\ &\qquad \left. + f'(y_u) A_{u,u+\varepsilon 2^{-n}} \right] \mathrm{d}u \\ &= \sum_{k=0}^{[1/\varepsilon]2^n - 1} \int_0^1 \left[ \frac{f(y_k^n) + f(y_{k+1}^n)}{2} (x_{k+1}^n - x_k^n) + f'(y_k^n) A_{k,k+1}^n \right] \mathrm{d}u, \end{aligned}$$

where we have written $x_k^n = x_{\varepsilon 2^{-n}(k+u)}$, $y_k^n = y_{\varepsilon 2^{-n}(k+u)}$ and $A_{k,\ell}^n = A_{\varepsilon 2^{-n}(k+u),\varepsilon 2^{-n}(\ell+u)}$ for simplicity. Dividing again in two, we find

$$\begin{aligned} I_{n+1}(\varepsilon) = &\sum_{k=0}^{[1/\varepsilon]2^n-1} \int_0^1 \left[ \frac{f(y_{2k}^{n+1}) + f(y_{2k+1}^{n+1})}{2} (x_{2k+1}^{n+1} - x_{2k}^{n+1}) + f'(y_{2k}^{n+1}) A_{2k,2k+1}^{n+1} \right] \mathrm{d}u \\ &+ \sum_{k=0}^{[1/\varepsilon]2^n-1} \int_0^1 \left[ \frac{f(y_{2k+1}^{n+1}) + f(y_{2k+2}^{n+1})}{2} (x_{2k+2}^{n+1} - x_{2k+1}^{n+1}) + f'(y_{2k+1}^{n+1}) A_{2k+1,2k+2}^{n+1} \right] \mathrm{d}u. \end{aligned}$$

On the other hand, after another change of variable, we can rewrite

$$\begin{aligned} I_n(\varepsilon) = &\frac{1}{2} \sum_{k=0}^{[1/\varepsilon]2^n-1} \int_0^1 \left[ \frac{f(y_{2k}^{n+1}) + f(y_{2k+2}^{n+1})}{2} (x_{2k+2}^{n+1} - x_{2k}^{n+1}) + f'(y_{2k}^{n+1}) A_{2k,2k+2}^{n+1} \right] \mathrm{d}u \\ &+ \frac{1}{2} \sum_{k=0}^{[1/\varepsilon]2^n-1} \int_0^1 \left[ \frac{f(y_{2k+1}^{n+1}) + f(y_{2k+3}^{n+1})}{2} (x_{2k+3}^{n+1} - x_{2k+1}^{n+1}) + f'(y_{2k+1}^{n+1}) A_{2k+1,2k+3}^{n+1} \right] \mathrm{d}u. \end{aligned}$$

Writing $J_n(\varepsilon) = I_{n+1}(\varepsilon) - I_n(\varepsilon)$, yields

$$\begin{aligned} J_n(\varepsilon) = &\frac{1}{2} \sum_{k=0}^{[1/\varepsilon]2^n-1} \int_0^1 \Bigg[ f'(y_{2k}^{n+1}) A_{2k,2k+1}^{n+1} + f'(y_{2k+1}^{n+1}) A_{2k+1,2k+2}^{n+1} \\ &\qquad - f'(y_{2k}^{n+1}) A_{2k,2k+2}^{n+1} + \frac{f(y_{2k}^{n+1}) + f(y_{2k+1}^{n+1})}{2} (x_{2k+1}^{n+1} - x_{2k}^{n+1}) \\ &\qquad + \frac{f(y_{2k+1}^{n+1}) + f(y_{2k+2}^{n+1})}{2} (x_{2k+2}^{n+1} - x_{2k+1}^{n+1}) \\ &\qquad\qquad - \frac{f(y_{2k}^{n+1}) + f(y_{2k+2}^{n+1})}{2} (x_{2k+2}^{n+1} - x_{2k}^{n+1}) \Bigg] \mathrm{d}u \\ &+ \frac{1}{2} \sum_{k=0}^{[1/\varepsilon]2^n-1} \int_0^1 \Bigg[ f'(y_{2k}^{n+1}) A_{2k,2k+1}^{n+1} + f'(y_{2k+1}^{n+1}) A_{2k+1,2k+2}^{n+1} \\ &\qquad - f'(y_{2k+1}^{n+1}) A_{2k+1,2k+3}^{n+1} + \frac{f(y_{2k}^{n+1}) + f(y_{2k+1}^{n+1})}{2} (x_{2k+1}^{n+1} - x_{2k}^{n+1}) \\ &\qquad + \frac{f(y_{2k+1}^{n+1}) + f(y_{2k+2}^{n+1})}{2} (x_{2k+2}^{n+1} - x_{2k+1}^{n+1}) \\ &\qquad\qquad - \frac{f(y_{2k+1}^{n+1}) + f(y_{2k+3}^{n+1})}{2} (x_{2k+3}^{n+1} - x_{2k+1}^{n+1}) \Bigg] \mathrm{d}u \end{aligned}$$

$$
\begin{aligned}
&= \frac{1}{2}\sum_{k=0}^{[1/\varepsilon]2^n-1}\int_0^1 \Bigg[ f'(y_{2k}^{n+1})A_{2k,2k+1}^{n+1} + f'(y_{2k}^{n+1})A_{2k+1,2k+2}^{n+1} \\
&\qquad + f'(y_{2k}^{n+1})A_{2k+2,2k}^{n+1} + \frac{f(y_{2k}^{n+1})+f(y_{2k+1}^{n+1})}{2}(x_{2k+1}^{n+1}-x_{2k}^{n+1}) \\
&\qquad + \frac{f(y_{2k+1}^{n+1})+f(y_{2k+2}^{n+1})}{2}(x_{2k+2}^{n+1}-x_{2k+1}^{n+1}) \\
&\qquad\qquad - \frac{f(y_{2k}^{n+1})+f(y_{2k+2}^{n+1})}{2}(x_{2k+2}^{n+1}-x_{2k}^{n+1})\Bigg]\,\mathrm{d}u \\
&\quad + O((\varepsilon 2^{-n})^{3\alpha-1}) \\
&\quad + \frac{1}{2}\sum_{k=0}^{[1/\varepsilon]2^n-1}\int_0^1 \Bigg[ f'(y_{2k+1}^{n+1})A_{2k+1,2k+2}^{n+1} + f'(y_{2k+1}^{n+1})A_{2k+2,2k+3}^{n+1} \\
&\qquad - f'(y_{2k+1}^{n+1})A_{2k+1,2k+3}^{n+1} + \frac{f(y_{2k+2}^{n+1})+f(y_{2k+3}^{n+1})}{2}(x_{2k+3}^{n+1}-x_{2k+2}^{n+1}) \\
&\qquad + \frac{f(y_{2k+1}^{n+1})+f(y_{2k+2}^{n+1})}{2}(x_{2k+2}^{n+1}-x_{2k+1}^{n+1}) \\
&\qquad\qquad - \frac{f(y_{2k+1}^{n+1})+f(y_{2k+3}^{n+1})}{2}(x_{2k+3}^{n+1}-x_{2k+1}^{n+1})\Bigg]\,\mathrm{d}u \\
&\quad + O((\varepsilon 2^{-n})^{\alpha\wedge(3\alpha-1)}) \\
&= \frac{1}{2}\sum_{k=0}^{[1/\varepsilon]2^n-1}\int_0^1 \Bigg[ f'(y_{2k}^{n+1})\mathrm{area}(T_{\gamma_{2k+2}\gamma_{2k+1}\gamma_{2k}}) \\
&\qquad + \frac{f(y_{2k+2}^{n+1})-f(y_{2k+1}^{n+1})}{2}(x_{2k}^{n+1}-x_{2k+1}^{n+1}) \\
&\qquad - \frac{f(y_{2k}^{n+1})-f(y_{2k+1}^{n+1})}{2}(x_{2k+2}^{n+1}-x_{2k+1}^{n+1})\Bigg]\,\mathrm{d}u \\
&\quad + \frac{1}{2}\sum_{k=0}^{[1/\varepsilon]2^n-1}\int_0^1 \Bigg[ f'(y_{2k+1}^{n+1})\mathrm{area}(T_{\gamma_{2k+3}\gamma_{2k+2}\gamma_{2k+1}}) \\
&\qquad + \frac{f(y_{2k+3}^{n+1})-f(y_{2k+2}^{n+1})}{2}(x_{2k+1}^{n+1}-x_{2k+2}^{n+1}) \\
&\qquad - \frac{f(y_{2k+1}^{n+1})-f(y_{2k+2}^{n+1})}{2}(x_{2k+3}^{n+1}-x_{2k+2}^{n+1})\Bigg]\,\mathrm{d}u \\
&\quad + O((\varepsilon 2^{-n})^{\alpha\wedge(3\alpha-1)})
\end{aligned}
$$

$$
\begin{aligned}
&= \frac{1}{2} \sum_{k=0}^{[1/\varepsilon]2^n-1} \int_0^1 f'(y_{2k}^{n+1}) \Big[ \operatorname{area}(T_{\gamma_{2k+2}\gamma_{2k+1}\gamma_{2k}}) \\
&\qquad + \frac{y_{2k+2}^{n+1} - y_{2k+1}^{n+1}}{2} (x_{2k}^{n+1} - x_{2k+1}^{n+1}) \\
&\qquad - \frac{y_{2k}^{n+1} - y_{2k+1}^{n+1}}{2} (x_{2k+2}^{n+1} - x_{2k+1}^{n+1}) \Big] \, \mathrm{d}u \\
&\quad + \frac{1}{2} \sum_{k=0}^{[1/\varepsilon]2^n-1} \int_0^1 f'(y_{2k+1}^{n+1}) \Big[ \operatorname{area}(T_{\gamma_{2k+3}\gamma_{2k+2}\gamma_{2k+1}}) \\
&\qquad + \frac{y_{2k+3}^{n+1} - y_{2k+2}^{n+1}}{2} (x_{2k+1}^{n+1} - x_{2k+2}^{n+1}) \\
&\qquad - \frac{y_{2k+1}^{n+1} - y_{2k+2}^{n+1}}{2} (x_{2k+3}^{n+1} - x_{2k+2}^{n+1}) \Big] \, \mathrm{d}u \\
&\quad + O((\varepsilon 2^{-n})^{\alpha \wedge (3\alpha - 1)}),
\end{aligned}
$$

where $\operatorname{area}(T_{abc})$ stands for the oriented area of the triangle $T_{abc}$, and where the correction terms come from (2.7), the $\mathscr{C}^2$-regularity of $f$ and the fact that $x, y$ are $\alpha$-Hölder. Now given that

$$
\tfrac{1}{2}[(y_c - y_b)(x_a - x_b) - (y_a - y_b)(x_c - x_b)] = \operatorname{area}(T_{abc}), \tag{2.11}
$$

we finally obtain

$$
I_{n+1}(\varepsilon) - I_n(\varepsilon) = O((\varepsilon 2^{-n})^{\alpha \wedge (3\alpha - 1)}),
$$

which yields the desired uniform convergence of $\{I_n, n \in \mathbb{N}\}$ toward some $I_\infty$. In addition, because $I_0(\varepsilon) = I_\varepsilon^\gamma(f)$, we have

$$
I_\infty(\varepsilon) = I_\varepsilon^\gamma(f) + O(\varepsilon^{\alpha \wedge (3\alpha - 1)}).
$$

This completes the proof in the case $m = 1$. Let us explain briefly the extension to the general case $m \geq 2$. Let $\Delta_n$ be the set of dyadics of order $n$ on $[0,1]$ and use the notation $t' = t + 2^{-n}$ and $\tau = \frac{t+t'}{2}$ for $t \in \Delta_n$. Let $\{w_n\}$ be the sequence defined by

$$
\begin{aligned}
w_n = {} & \sum_{t \in \Delta_n} (x_{t'} - x_t) \int_0^1 f((1-\alpha) y_t + \alpha y_{t'}) \nu_m(\mathrm{d}\alpha) \\
& + \sum_{k=0}^{2m-2} \frac{1}{(k+1)!} \sum_{t \in \Delta_n} f^{(k+1)}(y_t) A_{tt'}[(y - y_t)^k].
\end{aligned}
$$

Using a Taylor expansion, we can show that there exists a decomposition $w_{n+1} - w_n = U_n + V_n$ with $|U_n| \leq C2^{n(1-(2m+1)\alpha)}$ for some constant $C$ and

$$V_n = \sum_{k=0}^{2m-2} \frac{1}{(k+1)!}\left(\sum_{t\in\Delta_n}\{f^{(k+1)}(y_\tau)A_{\tau t'}[(y-y_\tau)^k] - f^{(k+1)}(y_t)A_{tt'}[(y-y_t)^k]\}\right).$$

Hence, $|V_n| \leq C2^{n(1-(2m+1)\alpha)}$ and the sequence $\{w_n\}$ converges absolutely. One can then finish the proof exactly as in the case $m = 1$. □

**Lemma 2.7.** *The function $I_\infty$ given in Lemma 2.6 is constant on* $[0, 1]$.

**Proof.** As in the proof of Lemma 2.6, we consider only the case $m = 1$ because the general case $m \geq 2$ can be handled completely analogously, with heavier notation. Once again, we set $A_{st}$ for $A_{st}(1)$. It is clear from the definition of $I_n$ and the uniqueness of the limit $I_\infty$ that

$$I_\infty(1) = I_\infty(2^{-1}) = I_\infty(2^{-2}) = \cdots = I_\infty(2^{-n}) = \cdots \tag{2.12}$$

for all $n \in \mathbb{N}$. We next prove that $I_\infty$ is constant on dyadics. From (2.12) and an induction argument, it suffices to prove that if $k2^{-n}$ and $(k+1)2^{-n}$ are two dyadics such that $I_\infty(k2^{-n}) = I_\infty((k+1)2^{-n}) = \ell$, then $I_\infty((k+1/2)2^{-n}) = \ell$. Using the notation $k_n^m = k2^{-(n+m)}$ we have, for all $m \in \mathbb{N}$,

$$\begin{aligned}
I_m\left(\left(k+\frac{1}{2}\right)2^{-n}\right) &= \frac{2^{n+m}}{2k+1}\int_0^1 [f(y_u) + f(y_{u+(k+1/2)_n^m})](x_{u+(k+1/2)_n^m} - x_u)\,\mathrm{d}u \\
&\quad + \frac{2^{n+m+1}}{2k+1}\int_0^1 f'(y_u)A_{u,u+(k+1/2)_n^m}\,\mathrm{d}u \\
&= \frac{2^{n+m}}{2k+1}\int_0^1 [f(y_{u+1_{n+1}^m}) + f(y_{u+(k+1)_n^m})](x_{u+(k+1)_n^m} - x_{u+1_{n+1}^m})\,\mathrm{d}u \\
&\quad + \frac{2^{n+m}}{2k+1}\int_0^1 f'(y_{u+1_{n+1}^m})A_{u+1_{n+1}^m,u+(k+1)_n^m}\,\mathrm{d}u + O(2^{-m\alpha}) \\
&= \frac{2k+2}{2k+1}I_m((k+1)2^{-n}) - \frac{1}{2k+1}I_{n+m+1}(1) + O(2^{-m[(3\alpha-1)\wedge\alpha]}) \\
&\quad + \frac{2^{n+m}}{2k+1}\int_0^1 [f(y_{u+1_{n+1}^m}) - f(y_u)](x_{u+(k+1)_n^m} - x_{u+1_{n+1}^m})\,\mathrm{d}u \\
&\quad - \frac{2^{n+m}}{2k+1}\int_0^1 [f(y_{u+(k+1)_n^m}) - f(y_{u+1_{n+1}^m})](x_{u+1_{n+1}^m} - x_u)\,\mathrm{d}u \\
&\quad + \frac{2^{n+m+1}}{2k+1}\int_0^1 f'(y_u)(A_{u+1_{n+1}^m,u+(k+1)_n^m} - A_{u,u+(k+1)_n^m} + A_{u,u+1_{n+1}^m})\,\mathrm{d}u
\end{aligned}$$

$$= \frac{2k+2}{2k+1} I_m((k+1)2^{-n}) - \frac{1}{2k+1} I_{n+m+1}(1) + O(2^{-m[(3\alpha-1)\wedge\alpha]}),$$

where the last line comes from (2.11). Making $m \to \infty$ yields

$$I_\infty\left(\left(k+\frac{1}{2}\right)2^{-n}\right) = \frac{2k+2}{2k+1}\ell - \frac{1}{2k+1}\ell = \ell,$$

which proves that $I_\infty$ is constant on the dyadics of $[0,1]$. Now because $I_n(\varepsilon)$ is obviously continuous in $\varepsilon$ and because the convergence in Lemma 2.6 is uniform, Dini's lemma entails that $I_\infty(\varepsilon)$ is continuous. Hence, $I_\infty$ is constant on $[0,1]$, as desired. □

**Proof of Theorem 2.5.** From Lemma 2.7 and (2.10), we have

$$I_\varepsilon^\gamma(f) = I_\infty(0) + O(\varepsilon^{((2m+1)\alpha-1)\wedge\alpha})$$

which, because $\alpha > 1/(2m+1)$, proves the convergence of $\{I_\varepsilon^\gamma(f), \varepsilon > 0\}$ toward some limit $I^\gamma$. Finally, because any $\alpha$-Lévy area of order $2n-2$ is also an $\alpha$-Lévy area of order $2m-2$ for $n > m$, the compatibility relationships follow straightforwardly from the fact that the correction term is in $\mathrm{O}(\varepsilon^{((2m+1)\alpha-1)\wedge\alpha})$. □

## 3. Differential equations in the case $\alpha > 1/3$

Recent works study equations of type (1.1) in the Russo–Vallois setting and in a Stratonovich sense. For example, in Errami and Russo [8], existence and uniqueness are proved for $H > 1/3$ with the following definition: a solution $X$ to (1.1) is a process such that $(X, B^H)$ is a symmetric vector cubic variation process and such that for every smooth $\varphi: \mathbb{R}^2 \to \mathbb{R}$ and every $t \geq 0$,

$$\int_0^t Z_s \,\mathrm{d}^\circ X_s = \int_0^t Z_s b(X_s)\,\mathrm{d}s + \int_0^t Z_s \sigma(X_s)\,\mathrm{d}^\circ B_s^H - \tfrac{1}{4}\int_0^t \sigma\sigma'(X_s)\,\mathrm{d}[Z, B^H, B^H]_s,$$

where $Z_s = \varphi(X_s, B_s^H)$ and $[Z, B^H, B^H]_s$ denotes the cubic covariation. In Nourdin [14], another type of equation is proposed, relying on the Newton–Côtes integrator and allowing any value of $H$, but the solution is supposed to be a priori of the type $X_s = h(B_s^H, V_s)$ with $V$ of bounded variation.

In this section we present yet another approach, which is more general and, hopefully, simpler. Fix $\alpha > \beta > 1/3$, a time horizon $T = 1$ and $x: [0,1] \to \mathbb{R} \in \mathrm{C}^\alpha$ once and for all. From Theorem 2.5, we know that

$$\int_0^1 f(y_s)\,\mathrm{d}^{A,1} x_s \stackrel{\text{def}}{=} \lim_{\varepsilon\to 0} \varepsilon^{-1} \int_0^1 \frac{f(y_s) + f(y_{s+\varepsilon})}{2}(x_{s+\varepsilon} - x_s)\,\mathrm{d}s + \varepsilon^{-1}\int_0^1 f'(y_s) A_{s,s+\varepsilon}\,\mathrm{d}s$$

exists as soon as $y: [0,1] \to \mathbb{R} \in \mathrm{C}^\beta$, $f: \mathbb{R} \to \mathbb{R} \in \mathscr{C}^2$ and $A$ is a $\beta$-Lévy area of order 0 associated with $\gamma = (x, y)$; that is, it verifies (2.8) and (2.9) with $\beta$ instead of $\alpha$.

Consider the formal equation

$$\mathrm{d}y_t = b(y_t)\,\mathrm{d}t + \sigma(y_t)\,\mathrm{d}x_t, \qquad t\in[0,1], \qquad y(0)=y_0\in\mathbb{R}, \tag{3.13}$$

where the unknown function $y$ is assumed to belong to $\mathrm{C}^\beta$. The following definition is inspired by rough paths theory:

***Definition 3.1.*** *A solution to (3.13) is a couple $(y,A)$ that verifies the following statements:*

- *$y:[0,1]\to\mathbb{R}$ belongs to $\mathrm{C}^\beta$.*
- *$A$ is a $\beta$-Lévy area of order 0 associated with $\gamma=(x,y)$.*
- *For all $t\in[0,1]$,*

$$y_t = y_0 + \int_0^t b(y_s)\,\mathrm{d}s + \int_0^t \sigma(y_s)\,\mathrm{d}^{A,1}x_s.$$

In this definition, we see that the meaning given to

$$\int_0^t \sigma(y_s)\,\mathrm{d}x_s$$

in (3.13) is contained in the concept of solution. We now state the main result of this section:

**Theorem 3.2.** *Let $\sigma:\mathbb{R}\to\mathbb{R}$ be a $\mathscr{C}^2$-function and let $b:\mathbb{R}\to\mathbb{R}$ be a Lipschitz function. Assume moreover that $\sigma,\sigma',\sigma''$ and $b$ are bounded. Then, for any $\beta\in(\frac{1}{3},\alpha)$, the equation (3.13) admits a unique solution $(y,A)$ in the sense of Definition 3.1.*

For the proof of Theorem 3.2 we need the following key lemma, which estimates the behavior of $\int_s^t f(y_u)\,\mathrm{d}^{A,1}x_u$ when $t-s$ is small:

**Lemma 3.3.** *Let $\beta\in(\frac{1}{3},\alpha)$, $y:[0,1]\to\mathbb{R}\in\mathrm{C}^\beta$, $A$ be a $\beta$-Lévy area of order 0 associated with $\gamma=(x,y)$ and $f:\mathbb{R}\to\mathbb{R}$ be a $\mathscr{C}^2$-function. Assume that $f,f'$ and $f''$ are bounded. Then, for all $s,t\in[0,1]$,*

$$\begin{aligned}&\left|\int_s^t f(y_u)\,\mathrm{d}^{A,1}x_u - \frac{f(y_s)}{t-s}\int_s^t (x_{u+t-s}-x_u)\,\mathrm{d}u - \frac{f'(y_s)}{t-s}\int_s^t A_{u,u+t-s}\,\mathrm{d}u\right|\\ &\quad\le c_{f,x}(|y|_\beta + |y|_\beta^2 + |y|_\beta|A|_{2\beta})(t-s)^{\alpha+\beta},\end{aligned} \tag{3.14}$$

*where $c_{f,x}$ denotes a constant that depends only on $f$ (and its derivatives) and $x$.*

**Proof.** Considering $\int_s^t$ instead of $\int_0^1$ in the proof of Lemma 2.6 (in the case $m=1$), we can write

$$\left|\int_s^t f(y_u)\,\mathrm{d}^{A,1}x_u - \frac{1}{t-s}\int_s^t \frac{f(y_u)+f(y_{u+t-s})}{2}(x_{u+t-s}-x_u)\,\mathrm{d}u\right.$$

$$-\frac{1}{t-s}\int_s^t f'(y_u)A_{u,u+t-s}\,\mathrm{d}u\bigg|$$
$$=|I_\infty(t-s)-I_0(t-s)|\le c_{f,x}(|y|_\beta|A|_{2\beta}+|y|_\beta^2)(t-s)^{3\beta},$$

but

$$\left|\frac{1}{t-s}\int_s^t\left(\frac{f(y_u)+f(y_{u+t-s})}{2}-f(y_s)\right)(x_{u+t-s}-x_u)\,\mathrm{d}u\right|$$
$$\le\frac{1+2^\beta}{2}|f'|_\infty|x|_\alpha|y|_\beta(t-s)^{\alpha+\beta}$$

and

$$\left|\frac{1}{t-s}\int_s^t(f'(y_u)-f'(y_s))A_{u,u+t-s}\,\mathrm{d}u\right|\le|f''|_\infty|y|_\beta|A|_{2\beta}|t-s|^{3\beta},$$

which entails (3.14). □

**Proof of Theorem 3.2.** Without loss of generality we can assume that $y_0=0$. Consider $E^\beta$ the set of couples $(y,A)$ with $y:[0,1]\to\mathbb{R}$ in $\mathrm{C}^\beta$ and $A$ a $\beta$-Lévy area of order 0 associated with $(x,y)$, endowed with the norm

$$N(y,A)=|y|_\beta+|A|_{2\beta}=\sup_{t\ne s}\frac{|y_t-y_s|}{|t-s|^\beta}+\sup_{t\ne s}\frac{|A_{st}|}{|t-s|^{2\beta}}<+\infty.$$

With this norm, $E^\beta$ is a Banach space. Besides, for every $\delta>0$, if $E_\delta^\beta$ denotes the set of restrictions of $(y,A)\in E^\beta$ to $[0,\delta]$, then $E_\delta^\beta$ is also a Banach space endowed with the norm $N$. Considering $(y,A)\in E_\delta^\beta$, define

$$\tilde{y}_t=\int_0^t\sigma(y_s)\,\mathrm{d}^{A,1}x_s+\int_0^t b(y_s)\,\mathrm{d}s,\qquad t\in[0,\delta],$$

and

$$\tilde{A}_{st}=\int_s^t x_u\sigma(y_u)\,\mathrm{d}^{A,1}x_u+\int_s^t x_u b(y_u)\,\mathrm{d}u-\tfrac{1}{2}(x_t+x_s)(\tilde{y}_t-\tilde{y}_s),\qquad(s,t)\in[0,\delta]^2.$$

Note here that, in the definition of $\tilde{A}$, we need to give a meaning to

$$\int_s^t x_u\sigma(y_u)\,\mathrm{d}^{A,1}x_u$$

and not only to $\int_s^t\sigma(y_u)\,\mathrm{d}^{A,1}x_u$. Following exactly the same arguments as in the proof of Theorem 2.5, there is no difficulty to consider integrals of this type and also to obtain the equivalent to Lemma 3.3 for this case.

Thanks to Lemma 3.3, we have

$$|\tilde{y}_t - \tilde{y}_s| \leq |\sigma|_\infty |x|_\alpha |t-s|^\alpha + |\sigma|_\infty |A|_{2\alpha} |t-s|^{2\alpha} + c_{\sigma,x}(|y|_\beta + |y|_\beta^2 + |y|_\beta |A|_{2\beta})(t-s)^{\alpha+\beta} + |b|_\infty |t-s|,$$

so that $\tilde{y} \in \mathrm{C}^\beta$. Condition (2.8) for $\tilde{A}$ is easily verified, while condition (2.9) with $\alpha$ replaced by $\beta$ also holds, because, using the equivalent of Lemma 3.3 for $\int_s^t x_u \sigma(y_u)\,\mathrm{d}^{A,1} x_u$, we can write

$$|\tilde{A}_{st}| \leq \tfrac{1}{2}|x|_\alpha^2 |\sigma|_\infty |t-s|^{2\alpha} + \tfrac{1}{2}|x|_\alpha |\sigma'|_\infty |A|_{2\beta} |t-s|^{2\beta+\alpha} + 2|x|_\infty |b|_\infty |t-s| + c_{\sigma,x}(|y|_\beta + |y|_\beta^2 + |y|_\beta |A|_{2\beta})(t-s)^{3\beta}.$$

In other words, $(\tilde{y}, \tilde{A}) \in E_\delta^\beta$ and we get

$$N(\tilde{y}, \tilde{A}) \leq c_{\sigma,b,x}(1 + \delta^{\alpha-\beta} N(y,A) + \delta^{\alpha-\beta} N(y,A)^2) \tag{3.15}$$

for a certain constant $c_{\sigma,b,x}$ that depends only on $\sigma$, $b$ and $x$. Because the set

$$\mathscr{U} = \{u \in \mathbb{R}_+^* : c_{\sigma,b,x}(1 + \delta^{\alpha-\beta} u + \delta^{\alpha-\beta} u^2) < u\}$$

is not empty as soon as $\beta$ is small enough, we deduce that the ball $E_\delta^\beta(R) = \{(y,A) \in E_\delta^\beta : N(y,A) \leq R\}$ is invariant by $T : E^\beta \to E^\beta$ defined by $T(y,A) = (\tilde{y}, \tilde{A})$ for $\beta$ small enough and $R \in \mathscr{U}$. Moreover, $T$ is actually a contraction provided $\delta$ is chosen small enough.

Given that we are working in $E_\delta^\beta(R)$, we can apply a standard fixed point argument for $T$, whose details are left to the reader. This leads to a unique solution to equation (3.13) on a small interval $[0, \delta]$. Notice now that the estimate (3.15) does not depend on the initial condition $y_0$, because $b$ and $\sigma$, together with its derivatives, are bounded. This makes it possible to obtain the unique solution on an arbitrary interval $[0, k\delta]$, $k \geq 1$, using a constant step $\delta$ and patching solutions on $[j\delta, (j+1)\delta]$. □

We now show how our concept of corrected symmetric integral extends the Russo–Vallois (RV) symmetric integral when the class of integrands is more specific. First, we define the Lévy area adapted to RV integrals:

**Lemma 3.4.** *Let $x : [0,1] \to \mathbb{R}$ be a function in $\mathrm{C}^\alpha$, $h : \mathbb{R}^2 \to \mathbb{R}$ be a $\mathscr{C}^{2,1}$-function and $\ell : [0,1] \to \mathbb{R}$ be a function of bounded variation. Define $y : [0,1] \to \mathbb{R}$ by $y_t = h(x_t, \ell_t)$. Then $y \in \mathrm{C}^\alpha$ and the Russo–Vallois symmetric integral $\int_r^s y\,\mathrm{d}^\circ x$ exists for all $r, s \in [0,1]$. Moreover, the function $A$ defined by*

$$A_{rs} = \int_r^s y\,\mathrm{d}^\circ x - \frac{y_r + y_s}{2}(x_s - x_r) \tag{3.16}$$

*is an $\alpha$-Lévy area of order 0 associated with $\gamma = (x, y)$, satisfying*

$$|A_{rs}| \le L|s - r|^{3\alpha} \tag{3.17}$$

*for all $r, s \in [0,1]$ and some universal constant $L$.*

**Proof.** For simplicity we consider only the case $y_t = h(x_t)$ with $h : \mathbb{R} \to \mathbb{R}$ a $\mathscr{C}^2$-function. The general case can be proved analogously. The fact that $y \in \mathrm{C}^\alpha$ and that $\int_r^s h(x)\,\mathrm{d}^\circ x$ exists for all $r, s \in [0,1]$ is well known and easy to obtain using a Taylor-type expansion. In addition, we know that $\int_r^s h(x)\,\mathrm{d}^\circ x = H(x_s) - H(x_r)$ for any primitive function $H$ of $h$. With the help of another Taylor expansion, it is then easy to show (3.17). Finally, the condition (2.6) is ascertained readily using the identity (2.11), which proves that $A$ is an $\alpha$-Lévy area of order 0 and finishes the proof of the lemma. □

The desired extension now follows readily:

**Proposition 3.5.** *With the same notation as in Lemma 3.4, we have*

$$\int_a^b f(y_s)\,\mathrm{d}^{A,1} x_s = \int_a^b f(y_s)\,\mathrm{d}^\circ x_s$$

*for any function $f : \mathbb{R} \to \mathbb{R}$ of class $\mathscr{C}^2$.*

**Proof.** Thanks to (3.17), we have

$$\lim_{\varepsilon \to 0} \varepsilon^{-1} \int_0^1 f'(y_u) A_{u,u+\varepsilon}\,\mathrm{d}u = 0,$$

which entails the required identification. □

***Remark 3.6.*** We do not know if it is possible to construct a Lévy area $A$ such that

$$\int_a^b f(y_s)\,\mathrm{d}^{A,m} x_s = \int_a^b f(y_s)\,\mathrm{d}^{NC,m} x_s$$

with the notation of Lemma 3.4, for any function $f : \mathbb{R} \to \mathbb{R}$ of class $\mathscr{C}^{2m}$, in the case $m \ge 2$. An area like

$$A_{rs}(y^q) = \frac{1}{q+1}\left(\int_r^s y_u^{q+1}\,\mathrm{d}^\circ x_u - (x_s - x_r)\int_0^1 (y_r + \theta(y_s - y_r))^{q+1} \nu_m(\mathrm{d}\theta)\right)$$

for $q \le m - 1$ would be the most natural candidate, but unfortunately only

$$|A_{st}[(y - \xi)^q]| \le c|t - s|^{3\alpha}$$

is fulfilled, in general, and not (2.7).

Finally, the next corollary shows that our solution process in Theorem 3.2 coincides with those given in Errami and Russo [8] or Nourdin [14] through a Doss–Sussmann's representation. If we could give a positive answer to the above remark, then the identification with Nourdin [14] would hold for all $m \geq 2$.

**Corollary 3.7.** *When $m=1$ and $\alpha > \beta > 1/3$, the unique solution $(y, A)$ to (3.13) can be represented as follows: The function $y:[0,1] \to \mathbb{R}$ is given by $y_t = u(x_t, a_t)$, where $u:\mathbb{R}^2 \to \mathbb{R}$ is the unique solution to*

$$\frac{\partial u}{\partial x}(x,v) = \sigma(u(x,v)) \quad \text{and} \quad u(0,v) = v \qquad \text{for all } v \in \mathbb{R}, \tag{3.18}$$

*and $a:[0,1] \to \mathbb{R}$ is the unique solution to*

$$\frac{\mathrm{d}a_t}{\mathrm{d}t} = \left\{ \frac{\partial u}{\partial a}(x_t, a_t) \right\}^{-1} b \circ u(x_t, a_t) \quad \text{and} \quad a_0 = y_0. \tag{3.19}$$

*The function $A$ is the $\beta$-Lévy area associated with $\gamma = (x, y)$ given by (3.16).*

**Proof.** It is clear that $y \in \mathrm{C}^\alpha \subset \mathrm{C}^\beta$ and we know from Proposition 3.5 that

$$\int_0^t \sigma(y_s)\, \mathrm{d}^{A,1} x_s = \int_0^t \sigma(y_s)\, \mathrm{d}^\circ x_s.$$

The easily established Itô–Stratonovich's formula verified by $\mathrm{d}^\circ$ shows that

$$u(x_t, a_t) = u(0, a_0) + \int_0^t \frac{\partial u}{\partial x}(x_s, a_s)\, \mathrm{d}^\circ x_s + \int_0^t \frac{\partial u}{\partial a}(x_s, a_s)\, \mathrm{d}a_s \tag{3.20}$$

for all $t \in [0,1]$. Hence, thanks to (3.18) and (3.19),

$$y_t = y_0 + \int_0^t \sigma(y_s)\, \mathrm{d}^\circ x_s + \int_0^t b(y_s)\, \mathrm{d}s = y_0 + \int_0^t b(y_s)\, \mathrm{d}s + \int_0^t \sigma(y_s)\, \mathrm{d}^{A,1} x_s$$

and, consequently, $(y, A)$ is the solution to (3.13). □

## Acknowledgements

We thank Denis Feyel for his helpful comments on the first version of this article, and two referees for very careful reading.

## References

[1] Alòs, E. and Nualart, D. (2002). Stochastic integration with respect to the fractional Brownian motion. *Stoch. Stoch. Rep.* **75** 129–152. MR1978896

[2] Cheridito, P. (2003). Arbitrage in fractional Brownian motion models. *Finance and Stochastics* **7** 533–553. MR2014249

[3] Cheridito, P. and Nualart, D. (2005). Stochastic integral of divergence type with respect to fractional Brownian motion with Hurst parameter $H \in (0, 1/2)$. *Ann. Inst. H. Poincaré Probab. Statist.* **41** 1049–1081. MR2172209

[4] Comte, F. and Renault, E. (1998). Long memory in continuous time volatility models. *Math. Finance* **8** 291–323. MR1645101

[5] Coutin, L. (2007). An introduction to (stochastic) calculus with respect to fractional Brownian motion. *Séminaire de Probabilités XL*. To appear.

[6] Coutin, L. and Qian, Z. (2002). Stochastic analysis, rough path analysis and fractional Brownian motions. *Probab. Theory Related Fields* **122** 108–140. MR1883719

[7] Cutland, N., Kopp, P. and Willinger, W. (1995). Stock price returns and the Joseph effect: a fractional version of the Black–Schole model. In *Seminar on Stochastic Analysis, Random Fields and Applications*, *Progr. Probab.* **36** 327–351. MR1360285

[8] Errami, M. and Russo, F. (2003). $n$-Covariation, generalized Dirichlet processes and calculus with respect to finite cubic variation processes. *Stochastic Process. Appl.* **104** 259–299. MR1961622

[9] Feyel, D. and De la Pradelle, A. (2006). Curvilinear integral along enriched paths. *Electron. J. Probab.* **11** 860–892. MR2261056

[10] Gradinaru, M., Nourdin, I., Russo, F. and Vallois, P. (2005). $m$-Order integrals and Itô's formula for non-semimartingale processes; the case of a fractional Brownian motion with any Hurst index. *Ann. Inst. H. Poincaré Probab. Statist.* **41** 781–806. MR2144234

[11] Lyons, T.J. (1998). Differential equations driven by rough signals. *Rev. Mat. Iberoamericana* **14** 215–310. MR1654527

[12] Lyons, T.J. and Qian, Z. (2003). *System Controls and Rough Paths.* Oxford: Oxford University Press. MR2036784

[13] Neuenkirch, A. and Nourdin, I. (2007). Exact rate of convergence of some approximation schemes associated with SDEs driven by a fBm. *J. Theoret. Probab.* To appear.

[14] Nourdin, I. (2007). A simple theory for the study of SDEs driven by a fractional Brownian motion, in dimension one. *Séminaire de Probabilités XLI*. To appear.

[15] Nourdin, I. and Simon, T. (2006). On the absolute continuity of one-dimensional SDE's driven by a fractional Brownian motion. *Statist. Probab. Letters* **76** 907–912. MR2268434

[16] Nualart, D. (1995). *The Malliavin Calculus and Related Topics.* New York: Springer. MR1344217

[17] Nualart, D. (2003). Stochastic calculus with respect to the fractional Brownian motion and applications. *Contemp. Math.* **336** 3–39. MR2037156

[18] Nualart, D. and Rasçanu, A. (2002). Differential equations driven by fractional Brownian motion. *Collect. Math.* **53** 55–81. MR1893308

[19] Russo, F. and Vallois, P. (1993). Forward, backward and symmetric stochastic integration. *Probab. Theory Related Fields* **97** 403–421. MR1245252